\documentclass[amssymb,amsfonts]{amsart}
\usepackage{amssymb, latexsym}

\newcommand{\half}{\frac{1}{2}}

\newcommand{\bp}{{\bf{p}}}

\newcommand{\be}{{\bf{e}}}

\def\be#1{\begin{equation} \label{#1}}
\def\bi{\begin{itemize}}
\def\bs{\begin{split}}
\def\im{{\text{Im}}}
\def\es{\end{split}}
\def\rr{{\mathbb{R}}}
\def\cc{{\mathbb{C}}}
\def\ba{\begin{align}}
\def\bas{\begin{align*}}
\def\philambda{\phi^{(\lambda)}}

\newcommand{\nabb}{\mbox{$\nabla \mkern-13mu /$\,}}

\def\ea{\end{align}}
\def\eas{\end{align*}}

\def\Re{{\hbox{Re}}}
\def\fatD{{\langle \nabla \rangle}}

\def\C{{\hbox{\mathbb C}}} 
\def\termone{{\text{Term}_{{1}} }}
\def\termtwo{{\text{Term}_{{2}} }}
\def\R{{\mathbb R}}

\def\C{{{\mathbb C}}} 

\def\textbf#1{{\bf #1}}

\newtheorem{theorem}{Theorem}
\theoremstyle{definition}

\theoremstyle{remark}

\theoremstyle{proposition}
\newtheorem{proposition}{Proposition}
\theoremstyle{lemma}
\newtheorem{lemma}{Lemma}
\theoremstyle{corollary}
\newtheorem{corollary}{Corollary}

\numberwithin{equation}{section}
\numberwithin{lemma}{section}
\numberwithin{remark}{section}
\numberwithin{theorem}{section}
\numberwithin{corollary}{section}
\numberwithin{proposition}{section}
\numberwithin{definition}{section}

\begin{document}

\title[Scattering for 3d NLS Below Energy]{Global existence and scattering  \\
 for rough solutions of a Nonlinear Schr\"odinger Equation on $\rr^3$}

\vspace{-0.3in}

\author{J. Colliander}
\thanks{J.C. is supported in part by N.S.F. grant DMS 0100595 and N.S.E.R.C. grant RGPIN 250233-03.}
\address{\small University of Toronto}

\author{M. Keel}
\thanks{M.K. was supported in part by the McKnight and Sloan Foundations.}
\address{\small University of Minnesota, Minneapolis}

\author{G. Staffilani}
\thanks{G.S. was supported in part by N.S.F. Grant
DMS 0100375 and the Sloan Foundation.}
\address{\small Massachusetts Institute of Technology}

\author{H. Takaoka}
\address{\small Kobe University}
\thanks{H.T. was supported in part by J.S.P.S. Grant No. 15740090}

\author{T. Tao}
\thanks{T.T. is a Clay Prize Fellow and was supported in part by
a grant from the Packard Foundation.}
\address{\small University of California, Los Angeles}


\begin{abstract}{We prove global existence and scattering for the defocusing, cubic nonlinear Schr\"odinger
equation in $H^s(\rr^3)$ for $s > \frac{4}{5}$.  The main new estimate in the argument is a Morawetz-type
inequality for the solution $\phi$.  This estimate bounds  $\|\phi(x,t)\|_{L^4_{x,t}(\rr^3 \times \rr)}$, whereas the well-known Morawetz-type estimate
of Lin-Strauss controls
$\int_0^{\infty}\int_{\rr^3}\frac{(\phi(x,t))^4}{|x|} dx dt $.
}
\end{abstract}

\maketitle


\section{Introduction and Statement of Results}

We study the following initial value problem for a cubic
defocusing nonlinear Schr\"odinger equation,
\begin{align}
\label{nls} i \partial_t \phi(x,t) + \Delta \phi (x,t)& =
|\phi(x,t)|^2 \phi(x,t), \quad
x \in \rr^3, t \geq 0,   \\
\phi(x,0) & = \phi_0(x) \; \in  H^s(\rr^3). \label{nlsdata}
\end{align}
Here $H^s(\rr^3)$ denotes the usual inhomogeneous Sobolev space.

It is known \cite{cw:local} that \eqref{nls}-\eqref{nlsdata} is
well-posed locally in time in $H^s(\rr^3)$ when\footnote{In
addition, there are local in time solutions from $H^{\frac{1}{2}}$
data, however, the time interval of existence depends upon the
profile of the initial data and not just upon the data's Sobolev
norm. Note that the  $\dot{H}^{\frac{1}{2}}(\rr^3)$ norm is {\em
critical} in the sense that it is invariant under the natural
scaling of solutions to \eqref{nls}.
} $s >\frac{1}{2}$ . In addition,  these local solutions enjoy
$L^2$ conservation,
\begin{align}
\label{l2conservation} ||\phi(\cdot, t)||_{L^2(\rr^3)} & =
||\phi_0(\cdot)||_{L^2(\rr^3)},
\end{align}
and the $H^1(\rr^3)$ solutions have the following conserved
energy,
\begin{align}
\label{energy} E(\phi)(t) & \equiv \int_{\rr^3} \frac{1}{2}
|\nabla_ x\phi(x,t)|^2 + \frac{1}{4} |\phi(x,t)|^4\ dx  =
E(\phi)(0).
 \end{align}

Together, these conservation laws and the local-in-time theory
immediately yield global-in-time well-posedness of
\eqref{nls}-\eqref{nlsdata} from data in $H^s(\rr^3)$ when $s \geq
1$. It is conjectured that \eqref{nls}-\eqref{nlsdata} is in fact
globally well-posed in time from all data included in the local
theory.  Previous work (\cite{firstnls}, extending
\cite{BScatter}) established this global theory when $s >
\frac{5}{6}$.
Our first goal here is to loosen  further the regularity
requirements on the initial data which ensure global-in-time
solutions. In addition we aim to loosen the symmetry assumptions
on the data which were previously used \cite{BScatter} to prove
scattering for rough solutions.

Before stating our main result, we recall some terminology (see
e.g. \cite{cazbook, GVScatter}). Write $S^L(t)$ for the flow map
$e^{it \Delta}$ corresponding to the linear Schr\"odinger
equation, and $S^{NL}(t)$ for the nonlinear flow, that is
$S^{NL}(t) \phi_0 = \phi(x,t)$ with $\phi, \phi_0$ as in
\eqref{nls},\eqref{nlsdata}. Given a solution\footnote{We can
easily extend the solution in \eqref{nls} to negative times by the
equation's time reversibility.} $\phi \in C\left((-\infty,
\infty), H^s(\rr^3)\right)$ of \eqref{nls}-\eqref{nlsdata}, define
the asymptotic states $\phi^{\pm}$  and  wave operators
$\Omega^{\pm}: H^s(\rr^3) \rightarrow H^s(\rr^3)$ by
\begin{align}
\phi^{\pm} & = \lim_{t \rightarrow \pm \infty} S^L(-t) S^{NL}(t)
\phi_0
\label{asymptoticstates}\\
\Omega^{\pm} \phi^{\pm} & = \phi_0  \label{waveoperators}
\end{align}
in so far as these limits exist in $H^s(\rr^3)$.  When the wave
operators $\Omega^{\pm}$ are surjective we say that
\eqref{nls}-\eqref{nlsdata} is {\em asymptotically complete} in
$H^s(\rr^3)$.

Our main result is the following:

\begin{theorem}  \label{maintheorem}
The initial value problem \eqref{nls}-\eqref{nlsdata} is
globally-well-posed from data $\phi_0 \in H^s(\rr^3)$ when $s >
\frac{4}{5}.$
 In addition, there is scattering
for these solutions.  More precisely, the wave operators
\eqref{waveoperators} exist and there is asymptotic completeness
on all of $H^s(\rr^3)$.
\end{theorem}

By {\em globally-well-posed}, we mean that given data $\phi_0 \in
H^s(\rr^n)$ as above, and any time $T > 0$,  there is a unique
solution to \eqref{nls}-\eqref{nlsdata}
\begin{align}
\phi(x,t) \in C([0,T]; H^s(\rr^n)) \label{summary}
\end{align}
which depends continuously in \eqref{summary} upon $\phi_0 \in
H^s(\rr^n)$.

We sketch the relationship of our results here with previous work.

Scattering in the space $H^{1}(\rr^3)$ was shown in
\cite{GVScatter}. Theorem \ref{maintheorem} extends part of the
work\footnote{In \cite{BScatter, bourgbook}, it is also shown that
the difference between the linear and nonlinear evolutions from
rough data has finite energy.  Our technique neither employs nor
implies such smoothing.} in \cite{BScatter, bourgbook}  where
global well-posedness was shown for general $H^s(\rr^3)$ data, $s
>\frac{11}{13}$.  (See \cite{bourg1} for a related result in two
space dimensions.)  In the case of radially symmetric data,
\cite{BScatter, bourgbook} establish global well-posedness and
scattering
for $\phi_0 \in H^s(\rr^3), s > \frac{5}{7}$. Theorem
\ref{maintheorem} also extends the result of \cite{firstnls},
where we showed global existence for $s > \frac{5}{6}$, with  no
scattering statement.


As in \cite{firstnls}, our arguments here preclude growth  of
$\|\phi(t)\|_{H^s(\rr^3)}$ by showing that the energy of a
smoothed version of the solution is almost conserved\footnote{The
phrase {\em{almost conserved}} is made precise in Proposition
\ref{almostconservation} below.}. We refer to  \cite{firstnls}
(pages 2-3) for remarks comparing the almost conservation law
approach used here with the argument in \cite{BScatter,
bourgbook}. See \cite{wavemaps,mkg, firstkdv, iteamIV} for further
applications of almost conservation laws; and \cite{iteamderiv,
iteamIII, iteamI} for instances where the inclusion of correction
terms in the almost conserved energy leads to sharp results.
Unlike our work in \cite{firstnls}, where
$\|\phi(t)\|_{H^s(\rr^3)}$ was bounded polynomially in time, we
ultimately obtain here a uniform bound.  The main new estimate
allowing such a uniform bound  is the Morawetz-type estimate
\eqref{spacetimeLfour} for the solution $u$ of any relatively
general defocusing nonlinear Schr\"odinger equation, see
\eqref{gnls} below.  Besides yielding the scattering results which
come along with such a uniform $L^4_{x,t}$ bound, this new
estimate is also the ingredient which pushes the allowed
regularity in Theorem \ref{maintheorem} below our previously
obtained $s > \frac{5}{6}$. We do not expect our results here to
be sharp. For example, we hope to extend Theorem \ref{maintheorem}
to allow lower values of $s$, using the correction terms mentioned
above and multilinear estimates (stemming from e.g.
\cite{collianderfancy, taolove}) to more tightly bound the
increment in the almost-conserved quantity.

Theorem \ref{maintheorem} above, like the referenced work on
global rough solutions for other dispersive equations, has a
number of motivations. We mention here three. First and most
obviously, we aim to better understand the global in time
evolution properties of known local-in-time solutions.  Second,
our results for rough solutions yield polynomial in time
bounds\footnote{In this paper, we in fact get a {\em uniform}
bound on the growth.} for the growth of some below-energy Sobolev
norms of {\em{smooth}} solutions. Such bounds give, for example, a
qualitative understanding of how the energy in a smooth solution
moves from high frequencies to low frequencies\footnote{If one has
a smooth solution with large but finite energy, the below-energy
Sobolev norms could presumably start relatively small and grow
large when the low frequencies of the solution grow in (for
example) $L^2$, while the high frequencies decrease in $L^2$.  A
polynomial bound on the rough norm's growth puts limits on this
movement of energy from high to low frequencies.}.  Third, we hope
that the techniques developed for these subcritical, rough initial
data problems can be used to address open problems for relatively
smooth solutions.  For an immediate example, our arguments below
give a new proof of the  finite energy scattering result of
\cite{GVScatter}. Also, the bounds we obtain on the global
Schr\"odinger admissible space-time norms of the solution depend
polynomially on the energy of the initial data, whereas previous
bounds were exponential. (See the remark in \cite{BScatter}, page
276, and \eqref{spacetimeLfour}, \eqref{almostclosed} below.)
There are of course more significant examples\footnote{Note added
in proof:  In the recent paper \cite{bigpaper},  we show global
well-posedness and scattering for the energy-critical (quintic)
defocusing analogue of \eqref{nls} from data in $H^s(\R^3), s \geq
1$.   The argument involves a frequency localized version of the
interaction Morawetz estimate (see Corollary \ref{herehere} below)
which holds for certain (hypothetical) blow-up solutions of the
quintic equation in three space dimensions. The argument also
relies on an almost conservation law for the frequency localized
mass of such solutions which is similar in spirit to Proposition
\ref{almostconservation} below.}  where low-regularity techniques
have helped to solve open problems for smooth solutions, e.g.
\cite{BRadial, TWave}.

The paper is organized as follows. In Section 2, after recalling
the standard Morawetz-type estimates from Lin-Strauss
\cite{linstrauss}, we introduce a Morawetz interaction potential
and prove it is bounded and monotone increasing. As a consequence,
we obtain the aforementioned spacetime $L^4_{x t}$ bound on
solutions of \eqref{nls}. Section 3 revisits the almost
conservation law argument in \cite{firstnls}, now in the setting
of an a-priori  $L^4_{x,t}$ bound on a spacetime slab. In Section
4, we first show in Proposition 4.1 how the almost conservation
law (Proposition 3.1), the interaction Morawetz inequality
\eqref{spacetimeLfour}, and the assumption $s > \frac{4}{5}$
combine with a scaling and bootstrap argument to give a
 uniform bound on $\|\phi(t)\|_{H^s(\rr^3)}$ and
the finiteness of $\|\phi\|_{L^4(\rr^3 \times [0, \infty))}$. The
scattering claims in Theorem \ref{maintheorem} follow from these
uniform bounds and by now well-known  arguments from earlier
scattering results of Brenner, Ginibre, Glassey, Morawetz,
Strauss,  and Velo (see surveys in \cite{cazbook,straussbook}).

Note that for finite energy solutions, that is $s = 1$,
Proposition 4.1 follows immediately from energy conservation and
the interaction Morawetz inequality \eqref{spacetimeLfour}.  Hence
in case $s = 1$, the arguments in Section 2 and the later part of
Section 4 below give a new, relatively direct proof of scattering
for \eqref{nls} in the energy class $H^1(\rr^3)$.  This result was
first established by Ginibre-Velo \cite{GVScatter}.

We conclude this introduction by setting some notation and
recalling the Strichartz estimates for the linear Schr\"odinger
operator on $\R^3$. Given  $A,B \geq 0$, we write $A \lesssim B$
to mean that for some universal constant $K > 2$, $A \leq K\cdot
B$.  We write $A \thicksim B$ when both $A \lesssim B$ and $B
\lesssim A$.  The notation $A \ll B$ denotes $B > K\cdot A$. We
write $\langle A \rangle \equiv (1 + A^2)^{\frac{1}{2}}$, and
$\fatD$ for the operator with Fourier multiplier $(1 +
|\xi|^2)^{\frac{1}{2}}$. The symbol $\nabla$ will denote the
spatial gradient. We will often use the notation $\frac{1}{2}+
\equiv \frac{1}{2} + \epsilon$ for some universal $0 < \epsilon
\ll 1$. Similarly, we write $\frac{1}{2}- \equiv \frac{1}{2} -
\epsilon$ .

Given Lebesgue space exponents $q,r$ and a function $F(x,t)$ on
$\rr^{n+1}$, we write
\begin{align}
\label{mixedlebesgue} ||F||_{L^q_tL^r_x(\rr^{n+1})} & \equiv
\left( \int_{\rr} \left( \int_{\rr^n} |F(x,t)|^r dx
\right)^{\frac{q}{r}} dt \right)^{\frac{1}{q}}.
\end{align}
This norm will be shortened to $L^q_tL^r_x$ for readability, or to
$L^r_{x,t}$ when $q=r$.

The Strichartz estimates involve the following definition:  a pair
of Lebesgue space exponents are called {\em{ Schr\"odinger
admissible}} for $\rr^{3+1}$ when $q,r \geq 2$, and
\begin{align}
\label{sa} \frac{1}{q} + \frac{3}{2r} & = \frac{3}{4}.
\end{align}

\begin{proposition}[Strichartz estimates in 3 space dimensions
(See e.g. \cite{segal, strichartz, GVStrichartz, yajima,
endpoint})] Suppose that $(q,r)$ and $(\tilde{q}, \tilde{r})$ are
{\em any} two Schr\"odinger admissible pairs as in \eqref{sa}.
Suppose too that $\phi(x,t)$ is a (weak) solution to the problem
\begin{align*}
(i \partial_t + \Delta) \phi(x,t) & = F(x,t), ~ (x,t) \in \rr^3 \times [0,T],  \\
\phi(x,0) & = \phi_0(x),
\end{align*}
for some data $u_0$ and $T > 0$.  Then  we have the estimate
\begin{align}
\label{strichartzestimate} ||\phi||_{L^q_tL^r_x([0,T] \times
\rr^3)} & \lesssim ||\phi_0||_{L^2(\rr^3)} +
||F||_{L^{\tilde{q}'}_t L^{\tilde{r}'}_x([0,T] \times \rr^3)}.
\end{align}
where $\frac{1}{\tilde{q}} + \frac{1}{\tilde{q}'} = 1,
\frac{1}{\tilde{r}} + \frac{1}{\tilde{r}'} =1$.
\end{proposition}

\section{The Morawetz interaction potential and a
spacetime $L^4$ estimate}

This section introduces an interaction potential generalization of
the classical Morawetz action and associated inequalities. We
first recall the standard Morawetz action centered at a point and
the proof that this action is monotonically increasing with time
when the nonlinearity is defocusing. The interaction
generalization is introduced in the second subsection. The key
consequence of the analysis in this section for the scattering
result is the $L^4_{x,t}$ estimate \eqref{spacetimeLfour}.

The discussion in this section will be carried out in the context
of the following  generalization of \eqref{nls}-\eqref{nlsdata}:
 \begin{align}
  \label{gnls}
i \partial_t u + \alpha \Delta u & = \mu f(|u|^2) u, &u:\R \times \R^3 \longmapsto \C, \\
u(0) & = u_0. \label{gnlsdata}
\end{align}
Here $f$ is a smooth function $f: \R^+ \longmapsto \R^+$ and
$\alpha$ and $\mu$ are real constants that permit us to easily
distinguish in the analysis below those terms arising from the
Laplacian or the nonlinearity. We also define $F(z) = \int_0^z
f(s) ds$.

We will use polar coordinates $x = r \omega,~ r > 0, \omega \in
S^2$, and write $\Delta_{\omega}$ for the Laplace-Beltrami
operator on $S^2$. For ease of reference below, we record some
alternate forms of the equation in \eqref{gnls}:
\begin{equation}
  \label{usubt}
   u_t = i \alpha \Delta u - i \mu f(|u|^2) u,
\end{equation}
\begin{equation}
  \label{ubarsubt}
  {\overline{u}}_t = - i \alpha \Delta {\overline{u}} + i \mu f(|u|^2) {\overline{u}},
\end{equation}
\begin{equation}
  \label{polargnls}
  u_t = i \alpha  u_{rr} + i \frac{2 \alpha }{r}  u_r
+ i \frac{\alpha}{r^2} \Delta_{\omega} u - i \mu f(|u|^2) u,
\end{equation}
\begin{equation}
  \label{rusubt}
   (ru_t) = i \alpha (r u )_{rr} + i \frac{\alpha}{r} \Delta_{\omega} u
- i \mu r f( |u|^2) u,
\end{equation}
\begin{equation}
  \label{rubarsubt}
   (r \overline{u}_t ) = - i \alpha (r {\overline{u}} )_{rr} - i \frac{\alpha}{r}
\Delta_{\omega} {\overline{u}} + i \mu f(|u|^2) \overline{u}.
\end{equation}

\subsection{Standard Morawetz action and inequalities}

We will call the following quantity the  {\it{Morawetz action
centered at 0}} for the solution $u$ of \eqref{gnls},
\begin{equation}
  \label{Mzero}
  M_0 [u] (t) = \int_{\R^3} \im [ {\overline{u}} (t,x) \nabla u(t,x) ] \cdot
\frac{x}{|x|} dx.
\end{equation}
We check using the equation that,
\begin{equation}
 \partial_t ( |u|^2) = - 2 \alpha \nabla \cdot   \im [ {\overline{u}} (t,x) \nabla u(t,x) ],
 \label{masscurrent}
 \end{equation}
 hence  we may interpret $M_0$ as the spatial average of the
radial component of the $L^2$-mass current. We might expect that
$M_0$ will increase with time if the wave $u$ scatters since such
behavior involves a broadening redistribution of the $L^2$-mass.
The following proposition of Lin and Strauss indeed gives
$\frac{d}{dt} M_0[u](t) \geq 0$ for defocusing equations.

\begin{proposition} \cite{linstrauss}
  If $u$ solves \eqref{gnls}-\eqref{gnlsdata} then the Morawetz action at 0 satisfies the
identity
\begin{equation}
  \label{Mzeroidentity}
  \partial_t M_0 [u] (t) = 4 \pi \alpha |u(t,0)|^2 + \int_{\R^3}
\frac{2 \alpha}{|x|} |\nabb_0 u(t,x) |^2 dx + \mu \int_{\R^3}
\frac{2}{|x|} \left\{ |u|^2 f(|u|^2 ) (t) - F(|u|^2 ) \right\} dx.
\end{equation}
where $\nabb_0$ is the angular component of the derivative,
\begin{equation}
  \label{nabbzero}
  \nabb_0 u = \nabla u - \frac{x}{|x|} (\frac{x}{|x|} \cdot \nabla u ).
\end{equation}
In particular, $M_0$ is an increasing function of time if the
equation \eqref{gnls} satisfies the repulsivity condition,
\begin{equation}
  \label{repulsivef}
  \mu \left\{ |u|^2 f(|u|^2 ) (t) - F(|u|^2 ) \right\} \geq 0.
\end{equation}
\end{proposition}

Note that for pure power potentials  $F(x) = \frac{2}{p+1}
x^{\frac{p+1}{2}}$, where the nonlinear term in \eqref{gnls} is
$|u|^{p-1} u$, the function $|u|^2 f(|u|^2) - F(|u|^2) =
\frac{p-1}{2} F(|u|^2).$ Hence condition \eqref{repulsivef} holds.

\begin{proof}
 Clearly, we may write
  \begin{align}
   \label{addoneoverr}
M_0 (t) & = \im \int_{\R^3} {\overline{u}} (t,x) (\partial_r +
\frac{1}{r} )
u(t,x) dx \\
\label{distribr}
  & = \im \int_0^\infty \int_{S^2} {\overline{ru}} ( ru )_r d \omega dr,
\end{align}
since we are working in three space dimensions. Integrating by
parts and using the equation \eqref{rusubt} gives,
\begin{eqnarray*}
\frac{d}{dt} {M_0} &=& \int_0^\infty \int_{S^2} {\overline{(ru)}}
(ru_t)_r
+ {\overline{(ru_t)}} (ru)_r d\omega dr \\
& = & - 2 \im \int_0^\infty \int_{S^2} {\overline{(ru)_r}} (ru_t) d\omega dr \\
 & = & - 2 \im \int_0^\infty \int_{S^2}  {\overline{(ru)_r}}
 \left\{ i \alpha (r u )_{rr} + i \frac{\alpha}{r} \Delta_{\omega} u
 - i \mu r f( |u|^2) u \right\} d\omega dr \\
&=& - 2 \alpha \Re \int_0^\infty \int_{S^2}  {\overline{(ru)_r}}
(ru)_{rr} \,
 d\omega dr - 2 \alpha \Re \int_0^\infty \int_{S^2} {\overline{(ru)_r}}
  \frac{1}{r}\Delta_\omega u  \,  d\omega dr
  \\
  & & \quad \quad + 2 \mu \Re \int_0^\infty \int_{S^2}
{\overline{(ru)_r}} r f( |u|^2) u  ~
  d\omega dr \\
& = & I + II + III.
\end{eqnarray*}
These three terms are analyzed separately and lead to the three
terms on the right side of \eqref{Mzeroidentity}.

\noindent{{\underline{Term I}}:} Since $\partial_r | (ru)_r |^2 =
2 \Re { \overline{(ru)_r}} (ru)_{rr},$ the $r$ integration in Term
I equals $ | (ru)_r |^2 |^\infty_0 = - |u(t,0)|^2$ which accounts
for the first term in \eqref{Mzeroidentity}.

\noindent{{\underline{Term II}}:}  Write $\Delta_\omega =
\nabla_\omega \cdot \nabla_\omega$ and integrate by parts to get,
\begin{equation*}
 II \; = \;   \alpha \Re \int_0^\infty \int_{S^2} \left[ \partial_r
  |\nabla_\omega u |^2 + \frac{2}{r} |\nabla_\omega u |^2 \right]
  d\omega dr.
\end{equation*}
Since $|\nabla_\omega u | \thicksim r |\nabla u |$, we know that
$|\nabla_\omega u |$ vanishes at the origin. Therefore, the first
term integrates to zero.  Finally, we can reexpress the remaining
term as claimed  in \eqref{Mzeroidentity} by inserting $r^2$ in
the numerator and denominator and then absorbing two factors of
$r$ using $\nabla_\omega u = r \nabb_0 u$.

\noindent{{\underline{Term III}}:} We expand the integrand using
the Leibniz rule to find $( \overline{u} + r {\overline{u}}_r ) r
f( |u|^2 ) u = r |u|^2 f(|u|^2) + r^2 f(|u|^2) u
{\overline{u}}_r.$ The first of these terms is purely real valued.
The real part of the second term may be reexpressed using $2 \Re
f(|u|^2) u {\overline{u}}_r = [F(|u|^2)]_r. $ Upon integrating
this last term by parts with respect to $r$, we obtain the third
expression in \eqref{Mzeroidentity}.

The remaining claim in the Proposition follows directly from
\eqref{Mzeroidentity}.
\end{proof}

We may center the above argument at any other point $y \in R^3$
with corresponding results.  Toward this end, define the
{\it{Morawetz action centered at $y$}} to be,
\begin{equation}
  \label{My}
  M_y[u] (t) = \int_{\R^3} \im [ {\overline{u} }(x)
  \nabla
u(x) ] \cdot \frac{x-y}{|x-y|} dx .
\end{equation}
We shall often drop the $u$ from this notation, as we did
previously in writing $M_0(t)$.
\begin{corollary} \label{cor2.1}
If $u$ solves \eqref{gnls} the Morawetz action at $y$ satisfies
the identity
\begin{equation}
  \label{Myidentity}
  {\frac{d}{dt}{M_y}} = 4 \pi \alpha |u(t,y)|^2 + \int_{\R^3} \frac{2
  \alpha}{|x-y|}
|\nabb_y u(t,x) |^2 dx + \int_{\R^3} \frac{2 \mu }{|x-y|} \left\{
|u|^2 f(|u|^2) - F(|u|^2)
  \right\}
dx,
\end{equation}
where $\nabb_y u \equiv \nabla u - \frac{x-y}{|x-y|} \left(
\frac{x-y}{|x-y|} \cdot \nabla u \right)$. In particular, $M_y$ is
an increasing function of time if the nonlinearity satisfies the
repulsivity condition \eqref{repulsivef}.
\end{corollary}

Corollary \ref{cor2.1} shows that a solution is, on average,
repulsed from any fixed point $y$ in the sense that $M_y[u](t)$ is
increasing with time.

For our scattering results, we'll need the following pointwise
bound for $M_y[u](t)$.

\begin{lemma}  Assume $u$ is a solution of \eqref{gnls} and $M_y[u](t)$ as in \eqref{My}.  Then,
  \begin{equation}
    \label{Myupper}
   | M_y (t) | \lesssim {{\| u(t) \|}_{\dot{H}^\half_x}^2} .
  \end{equation}
\end{lemma}

 \begin{proof} Without loss of generality we take $y = 0$.
 This is a refinement of the easy bound using Cauchy-Schwarz $| M_y (t)
| \lesssim {{\| u (t) \|}_{L^2_x}} {{\| \nabla u(t) \|}_{L^2_x}}.$
By duality
 \begin{align*}
 | \, \im \int_{\rr^3} \overline{u(x,t)} \partial_r u(x,t) dx \, | &
 \leq \| u \|_{\dot{H}^{\frac{1}{2}}(\rr^3)}
 \cdot \| \partial_r u \|_{\dot{H}^{-\frac{1}{2}}(\rr^3)}.
 \end{align*}
 It suffices to show $\|\partial_r u \|_{\dot{H}^{-\frac{1}{2}}(\rr^3)}
 \leq \| u \|_{\dot{H}^{-\frac{1}{2}}(\rr^3)}$. By duality and the definition
 $\partial_r \equiv \frac{x}{|x|} \cdot \nabla$, it remains to prove,
 \begin{align}
 \label{unitbegone}
 \| \frac{x}{|x|} f \|_{\dot{H}^{\frac{1}{2}}(\rr^3)} & \leq \|  f \|_{\dot{H}^{\frac{1}{2}}(\rr^3)},
 \end{align}
 for any $f$ for which the right hand side is finite. Inequality \eqref{unitbegone} follows from interpolating
 between the following two bounds,
 \begin{align*}
 \|\frac{x}{|x|} f\|_{L^2(\rr^3)} & \leq \| f\|_{L^2(\rr^3)} \\
 \|\frac{x}{|x|} f\|_{\dot{H}^1(\rr^3)} & \lesssim \| f\|_{H^1(\rr^3)}
 \end{align*}
 the first of which is trivial, the second of which follows from Hardy's inequality,
 \begin{align*}
 \| \nabla \left( \frac{x}{|x|} f\right) \|_{L^2} & \leq \| \frac{x}{|x|} \cdot \nabla f\|_{L^2} + \|\frac{1}{|x|} f \|_{L^2}\\
 & \lesssim \| \nabla f\|_{L^2}.
 \end{align*}
 \end{proof}
The well-known Morawetz-type  inequalities which have proven
useful in proving local decay or scattering for \eqref{gnls} arise
by integrating  the identity \eqref{Mzeroidentity} or
\eqref{Myidentity} in time.  For nonlinear Schr\"odinger
equations, this argument appears in the work of Lin and Strauss
\cite{linstrauss}, who cite as motivation earlier work on
Klein-Gordon equations by Morawetz \cite{morawetz}.

\begin{corollary}[Morawetz inequalities \cite{linstrauss}]  Suppose $u$ solves
\eqref{gnls}-\eqref{gnlsdata}.  Then for any $y \in \R^3$,
  \begin{equation}
    \label{MorawetzInequalities}
    2 \sup_{t \in [0,T]} {{\| u(t) \|}_{\dot{H}^\half_x}^2} \gtrsim
4 \pi \alpha \int_0^T |u(t,y)|^2 dt + \int_0^T \int_{\R^3} \frac{2
\alpha}{|x-y|} |\nabb_y u(t,x)|^2 dx dt
\end{equation}
\begin{equation*}
+ \int_0^T \int_{\R^3} \frac{2\mu}{|x-y|}  \left\{ |u|^2 f(|u|^2)
- F(|u|^2)
  \right\} dx dt.
  \end{equation*}
\end{corollary}

Assuming \eqref{gnls} has a repulsive nonlinearity as in
\eqref{repulsivef}, all terms on the right side of the inequality
\eqref{MorawetzInequalities} are positive.  The inequality
therefore gives in particular a bound uniform in $T$ for the
quantity $\int_0^T \int_{\R^3} \frac{|u(t,x)|^4}{|x-y|} dx dt$,
for solutions $u$ of \eqref{nls}.

In their proof of scattering in the energy space for the cubic
defocusing problem \eqref{nls}, Ginibre and Velo \cite{GVScatter}
combine this relatively localized\footnote{The bound mentioned
here may be considered localized since it implies decay of the
solution near the fixed point $y$, but doesn't preclude the
solution staying large at a point which moves rapidly away from
$y$, for example.} decay estimate with a bound surrogate for
finite propagation speed in order to show the solution is in
certain global-in-time Lebesgue spaces $L^q([0, \infty),
L^r(\rr^3))$.    Scattering follows rather quickly.

In the following section, we show how to establish an unweighted,
global in time Lebesgue space bound directly. The argument below
involves the identity \eqref{Myidentity}, but our estimate arises
eventually from the linear part of the equation, more specifically
from the first term on the right of \eqref{Myidentity}, rather
than the third (nonlinearity) term.

\subsection{Morawetz interaction potential}

Given a solution $u$ of \eqref{gnls}, we define the {\it{Morawetz
interaction potential}}  to be
\begin{equation}
  \label{M}
  M(t) = \int_{\rr^3} |u (t, y)|^2 M_y (t) dy.
\end{equation}
 The bound \eqref{Myupper} immediately implies
\begin{equation}
  \label{Mupper}
  |M(t)| \lesssim {{\| u (t) \|}_{L^2}^2}
 {{\| u(t) \|}_{\dot{H}^\half_x}^2} .
\end{equation}
If $u$ solves \eqref{gnls} then the identity \eqref{Myidentity}
gives us the following identity for $\frac{d}{dt} M(t)$,
  \begin{multline}
  \label{interactionsubt}
 \;\;\;  \frac{d}{dt}M(t) \; = \; 4 \pi \alpha \int_y |u(y)|^4 dy + \int_{\R^3} \int_{\R^3}
\frac{2 \alpha}{|x-y|} |u(y)|^2 |\nabb_y u(x) |^2 dx dy  \\
 \;\;\;\; + \int_{\R^3} \int_{\R^3} \frac{2 \mu }{|x-y|} |u(y)|^2 \left\{ |u(x)|^2 f(|u(x)|^2)
    - F(|u(x)|^2)   \right\} dx dy  \\ + \int_{\R^3} \partial_t ( |u(t,y)|^2 )  ~
M_y (t) dy.
\quad\quad\quad\quad\quad\quad\quad\quad\quad\quad\quad\quad\quad\quad\quad\quad\quad\quad\quad\quad\quad
\end{multline}
We write the right side of \eqref{interactionsubt} as $I + II +
III + IV$, and work now to rewrite this as a sum involving
nonnegative terms.

\begin{proposition}
\label{fournottoonegative}
  Referring to the terms comprising \eqref{interactionsubt}, we have
  \begin{equation}
   \label{fournottooneg}
   IV \geq - II.
  \end{equation}
Consequently, solutions of \eqref{gnls} satisfy
\begin{equation}
  \label{interactionmonotone}
  {\frac{d}{dt}{M(t)}} \geq 4 \pi \alpha \int_{\R^3} |u(t,y)|^4 dy + \int_{\R^3} \int_{\R^3}
  \frac{2\mu}{|x-y|}  |u(t ,y)|^2 \left\{ |u|^2 f(|u|^2)
    - F(|u|^2)   \right\} dx dy.
\end{equation}
In particular, $M(t)$ is monotone increasing for equations with
repulsive nonlinearities.
\end{proposition}

Assuming Proposition \ref{fournottooneg} for the moment, we
combine \eqref{Mupper} and \eqref{interactionmonotone} to obtain
the following estimate which plays the major new role in our
analysis in Sections 3 and 4 below,

\begin{corollary}  \label{herehere} Take $u$ to be a smooth solution to the initial value problem \eqref{gnls}-\eqref{gnlsdata}
above, under the repulsivity assumption \eqref{repulsivef}.  Then
we have the following {\it{interaction Morawetz inequalities}},
\begin{multline}
\label{intMorawineq} \;\;\;\;\;    2 {{\| u(0) \|}_{L^2}^2}
\sup_{t \in [0,T]} {{\| u(t) \|}_{\dot{H}_x^\half}^2}  \; \gtrsim
\; 4 \pi \alpha \int_0^T \int_{\rr^3}
|u(t,y)|^4 dy dt  \\
\;\; \; \; +  \int_0^T \int_y \int_x \frac{2 \mu }{|x-y|}
|u(t,y)|^2\left\{ |u|^2 f(|u|^2)
    - F(|u|^2)   \right\}(t,x) dx dy dt.
\end{multline}
In particular, we obtain the following spacetime $L^4([0, \infty)
\times \R^3)$ estimate,
\begin{equation}
  \label{spacetimeLfour}
  \int_0^T \int_{\R^3} |u(t,y)|^4 dy dt \lesssim {\| u_0 \|}_{L^2(\rr^3)}^2
\sup_{t \in [0,T]} {{\| u(t) \|}_{\dot{H}_x^\half}^2} .
\end{equation}
\end{corollary}
Of course, for solutions of \eqref{nls} starting from finite
energy initial data, the right side of \eqref{spacetimeLfour} is
uniformly bounded by energy considerations - leading to a rather
direct proof of the result in \cite{GVScatter} of scattering in
the energy space. This bound \eqref{spacetimeLfour} is also a key
part of our rough data scattering argument below.

\begin{proof} We now turn to the proof of Proposition \ref{fournottoonegative}. Use
  \eqref{masscurrent} to write
  \begin{eqnarray*}
    IV &=& - \int_{\R_y^3} \nabla \cdot \im [ 2 \alpha \overline{u}(y)
\nabla u(y) ]
M_y (t) dy \\
& = &- \int_y \int_x \partial_{y_l} \im [2 \alpha \overline{u}(y)
\partial_{y_l} u(y) ]~ \im [ \overline{u} (x) \frac{x_m - y_m}{|x-y|}
\partial_{x_m} u(x) ] dx dy,
  \end{eqnarray*}
 where repeated indices are implicitly summed. We integrate by parts in $y$, moving   the leading
 $\partial_{y_l}$ to the unit vector
 $\frac{x-y}{|x-y|}$. Note that,
  \begin{equation}
   \label{vectorid}
   \partial_{y_l} \left( \frac{x_m - y_m}{|x-y|} \right) =
   \frac{-\delta_{lm}}{|x-y|} + \frac{ (x_l - y_l)(x_m - y_m) }{|x-y|^3} .
 \end{equation}
Write  $\bp (x) = \im [ {\overline{u}} (x) \nabla u (x) ]$ for the
mass current at $x$ and use \eqref{vectorid} to obtain
\begin{equation}
  \label{projection}
  IV = -2 \alpha \int_y \int_x \left[  \bp (y) \cdot \bp (x) -
( \bp (y) \cdot \frac{x-y}{|x-y|} ) (\bp (x) \cdot
\frac{x-y}{|x-y|} ) \right] \frac{dx dy}{|x-y|}.
\end{equation}
The preceding integrand has a natural geometric interpretation. We
are removing the inner product of the components of $\bp (y)$ and
$\bp (x)$ parallel to the vector $\frac{x-y}{|x-y|}$ from the full
inner product of $\bp (y)$ and $\bp (x)$. This amounts to taking
the inner product of $\pi_{(x-y)^\perp} \bp (y) \cdot
\pi_{(x-y)^\perp} \bp (x)$ where we have introduced the
projections onto the subspace of $\R^3$ perpendicular to the
vector $\frac{x-y}{|x-y|}$. But
\begin{equation}
\label{crux} |\pi_{(x-y)^\perp} \bp (y)| = \; \big|\bp(y) -
\frac{x-y}{|x-y|} \big( \frac{x-y}{|x-y|} \cdot \bp(y) \big) \big|
\; = \; | \im [ \overline{u} (y ) \nabb_x u(y) | \leq |u(y)| \cdot
|\nabb_x u(y)|.
\end{equation}
A similar identity and inequality holds upon switching the roles
of $x$ and $y$ in \eqref{crux}. We have thus shown that
\begin{equation}
  \label{endingfour}
  IV \geq - 2 \alpha \int_y \int_x |u(x)| \cdot |\nabb_y u(x) | \cdot |u(y) |
 \cdot  |\nabb_x u(y)| \frac{dx dy}{|x-y|}.
\end{equation}
The conclusion \eqref{fournottooneg} follows by applying the
elementary bound
 $|a b | \leq \half (a^2 + b^2)$ with $a = |u(y)| \cdot |\nabb_y u(x)|$ and $ b = |u(x)|\cdot |\nabb_x u(y)|$.
\end{proof}



 \section{ Almost Conservation Law.}

 Keeping in mind that the energy  \eqref{energy} of our solutions might be infinite, our
 aim will be to control the growth in time of $E(I\phi)(t)$, where $I\phi$ is a smoothed
 version of $\phi$.  The operator $I$ depends on a parameter $N \gg 1$ to be chosen later, and the level
 of regularity $s <1$ at which we are working\footnote{We abuse notation and suppress this dependence,
 writing simply $I$ instead of $I_{s,N}$.}. We write,
 \begin{align}
 \label{Ioperator}
 \widehat{If}(\xi) & \equiv m_N(\xi) \hat{f}(\xi),
 \end{align}
 where
 the multiplier $m_N(\xi)$ is smooth, radially symmetric, nonincreasing
 in $|\xi|$ and
 \begin{align}
 \label{Iproperties}
 m_N(\xi) & =
 \begin{cases}
 1 & |\xi| \leq N \\
 \left( \frac{N}{|\xi|} \right)^{1 - s} & |\xi| \geq 2N.
 \end{cases}
 \end{align}
 The following two inequalities follow quickly from the definition of $I$, the $L^2$ conservation
 \eqref{l2conservation}, and by
 considering separately those frequencies $|\xi| \leq N$ and $|\xi| \geq N$.
 \begin{align} \label{boundenergy}
 E(I \phi)(t) & \lesssim \left(N^{1 - s}
 ||\phi(\cdot, t)||_{\dot H^s(\rr^3)} \right)^2 +
 ||\phi(t, \cdot)||_{L^4(\rr^3)}^4, \\
 \label{boundhsnorm}
 ||\phi(\cdot, t)||_{H^s(\rr^3)}^2 & \lesssim E(I \phi)(t) +
 ||\phi_0||_{L^2(\rr^3)}^2.
 \end{align}
 In studying the possible growth of our solution in time, we will estimate $E(I \phi)(t)$ rather than
 bounding $||\phi(t)||_{H^s(\rr^3)}$ directly.
 Of course, since \eqref{nls} is a nonlinear
 equation, $I\phi(x,t)$ is not a solution.  In particular, one
 doesn't expect $E(I\phi)(t)$ to be constant.  One of the main ingredients of Theorem \ref{maintheorem}
 is proving that this quantity is uniformly bounded in time. The local in time result which contributes to the
 proof of such a bound is what we mean
 by an {\em almost conservation law}.  Global well-posedness follows from \eqref{boundhsnorm}, a uniform
 bound on $E(I\phi)(t)$ in terms of $\|\phi_0\|_{H^s(\rr^3)}$, the fact that  \eqref{nls}-\eqref{nlsdata} is
 locally well posed when $s > \frac{1}{2}$,
 and a density argument.

\begin{proposition} [Almost Conservation Law] \label{almostconservation}
Assume we have  $s > \frac{1}{2},\; N \gg 1,$
$\phi_0 \in C^{\infty}_0(\rr^3)$, and a solution of
\eqref{nls}-\eqref{nlsdata} on a time  interval $[0,T]$  for which
\begin{align}
\label{smalll4onslab} ||\phi||_{L^4_{x,t}([0,T] \times \rr^3)} &
\lesssim \epsilon.
\end{align}
Assume in addition that $E(I\phi_0) \lesssim 1$.

We conclude  that for all $t \in [0,T]$,
\begin{align}
\label{increment} E(I \phi)(t) & =  E(I \phi_0) + O(N^{- 1+}).
\end{align}
\end{proposition}

\medskip

Equation \eqref{increment} asserts that $I \phi$, though not a
solution of the nonlinear problem \eqref{nls}, enjoys something
akin to energy conservation. If one could replace the increment
$N^{- 1+}$ in $E(I \phi)$ on the right side of \eqref{increment}
with $N^{- \alpha}$ for some $\alpha > 0$, one could repeat the
argument we give below to prove global well-posedness of
\eqref{nls}-\eqref{nlsdata} for all $s > \frac{3 + \alpha}{3 +
2\alpha}$. In particular, if $E(I \phi)(t)$ were conserved (i.e.
$\alpha = \infty$), one could show that
\eqref{nls}-\eqref{nlsdata} is globally well-posed when $s >
\frac{1}{2}$. Recall that the scale-invariant Sobolev space is
$\dot{H}^{\frac{1}{2}}(\rr^3)$.


Proposition \ref{almostconservation} is a modification of a
similar statement (also labelled Proposition 3.1) in
\cite{firstnls}. The statement in \cite{firstnls} establishes a
uniform time step, determined by the size of the modified energy
of the data $E(I \phi )$, on which there is almost conservation of
$E(I \phi )(t)$. Here we obtain an almost conservation property in
time intervals $[0,T]$  on which $\phi$ is assumed small in
$L^4_{x,t}$. Note that these intervals  may have various lengths,
and that the constant implicit in \eqref{increment} is independent
of these lengths.

The proof of Proposition \ref{almostconservation} proceeds by
pretending that $I \phi$ is a solution of $\eqref{nls}$ and using
the usual proof of energy conservation. We look at  the resulting
space-time integral in Fourier space, where we estimate various
frequency interactions separately.  In doing so, we'll need
control of a local-in-time norm $Z_I(t)$ involving the indices in
\eqref{sa},
\begin{align}
 \label{newz}
Z_I(t) & \equiv \sup_{q,r \; \hbox{admissible}} || \nabla I \phi
||_{L^q_tL^r_x([0,t] \times \rr^3)}
\end{align}   similar to those norms
that are usually bounded by the local in time existence theorem
for \eqref{nls}.  (See e.g. \cite{cw:local}). Since the norm here
includes the operator $I$, and as mentioned above, we will control
$Z_I(t)$ on time intervals of varying lengths, we think of the
following lemma as a modified local existence theory.



\begin{lemma} \label{interpolation}  Consider  $\phi(x,t)$ as in \eqref{nls}-\eqref{nlsdata}
defined on $[0, T^*] \times \rr^3$ where
\begin{align}
\label{gotitgoinon} \|\phi\|_{L^4_{x,t}([0,T^*] \times \rr^3)} &
\leq\epsilon,
\end{align}
for some universal constant $\epsilon$.  Assume too $\phi_0 \in
C^{\infty}_0(\rr^3)$.  Then for $s > \frac{1}{2}$ and sufficiently
large\footnote{Recall that $I \equiv I_{N,s}$ was defined in
\eqref{Ioperator}-\eqref{Iproperties}.} $N$,
\begin{align}
Z_I(T^*) & \leq C(||\phi_0||_{H^s(\rr^3)}). \label{boundnewz}
\end{align}
\end{lemma}

\noindent {\bf{Proof of Lemma \ref{interpolation}}}: Apply $I
\nabla$ to both sides of \eqref{nls}.  Choosing $\tilde{q}',
\tilde{r}' = \frac{10}{7}$, \eqref{strichartzestimate} and a
fractional Leibniz rule\footnote{Since $s > \frac{1}{2}$, the
multiplier for $\nabla^{\alpha} I$ is increasing in $|\xi|$ when
$\frac{1}{2} \leq \alpha \leq 1$.  Using this fact, one can easily
modify the usual proof of the fractional Leibniz rule so this rule
holds for the operators $\nabla^\alpha I$.  (See e.g. page 105 of
the exposition in \cite{taylorbook}, or  the articles \cite{cw},
\cite{katoponce}.)} give us that for all $0 \leq t \leq T$,
\begin{align*}
Z_I(t) & \lesssim ||\nabla I \phi_0||_{L^2(\rr^3)} + ||\nabla I
\phi ||_{L^{\frac{10}{3}}_{x,t}([0,t] \times \rr^3)} \cdot || \phi
||_{L^{5}_{x,t}([0,t] \times \rr^3)}^2.
\end{align*}
The $L^{\frac{10}{3}}$ factor here is bounded by $Z_I(t)$. We
claim that the remaining $L^5_{x,t}$ factors are bounded by,
\begin{align}
\label{interp} \|\phi\|_{L^5_{x,t}([0,T^*] \times \rr^3)} &
\lesssim \epsilon^{\delta_1} \cdot (Z_I(T^*))^{\delta_2}
\end{align}
for some $\delta_1, \delta_2 > 0$, and  $Z_I$ as in \eqref{newz}.
Assuming \eqref{interp} for the moment, we conclude that for $N$
sufficiently large,
\begin{align}
Z_I(t) & \lesssim 1 + \epsilon^{\delta_3} \left(Z_I(t)\right)^{1+
\delta_4}, \label{almostclosednew}
\end{align}
for some constants $\delta_3, \delta_4 > 0$.  For sufficiently
small choice of $\epsilon$, the bound \eqref{almostclosednew}
yields \eqref{boundnewz} for all $0 \leq t \leq T$, as desired.

It remains to prove \eqref{interp}. All space-time norms in this
proof will be taken on the slab $[0, T^*] \times R^3$, even when,
for legibility,  this isn't explicitly written. Write  $$\phi =
\psi_0 + \sum_{i=1}^\infty \psi_j$$ where $\psi_0$ has spatial
frequency support on $\langle \xi \rangle \lesssim N_1 \equiv N$
and the remaining $\psi_j$ each have dyadic spatial frequency
support $\langle\xi_j \rangle \sim N_j \equiv 2^{k_j},$ where $k_j
\gtrsim \log(N)$ are integers and  $j = 1,2, \ldots$. The argument
given below estimates the low frequency constituent $\psi_0$ with
the available $L^4$ and $L^{10}$ bounds; and the high frequency
pieces $\psi_j, j \geq 1$ with the $L^{\frac{10}{3}}$ and $L^{10}$
bounds.

Specifically, the definition of $I$ in \eqref{Iproperties} gives,
\begin{align*}
\|I\psi_j\|_{L^{10}_{x,t}} &\thicksim \begin{cases}
\|\psi_j\|_{L^{10}_{x,t}} & j = 0 \\
N^{1-s} (N_j)^{s-1} \|\psi_j\|_{L^{10}_{x,t}} & j = 1, 2, \ldots.
\end{cases}
\end{align*}
Using Sobolev's inequality, the left hand side here is bounded by
$Z_I(T^*)$.  Rewriting gives,
\begin{align}
\label{star} \|\psi_j\|_{L^{10}_{x,t}([0,T^*] \times \rr^3)} &
\lesssim \begin{cases}
Z_I(T^*) & j \, = \, 0 \\
N_j^{1-s} N^{s-1} Z_I(T^*) & j = 1,2,\ldots
\end{cases}.
\end{align}
Similarly,
\begin{align*}
\|\nabla I \psi_j \|_{L^{\frac{10}{3}}_{x,t}} & \thicksim N_j^s
N^{1-s} \|\psi_j \|_{L^{\frac{10}{3}}_{x,t}} \quad  j =
1,2,\ldots.
\end{align*}
Hence we get the following $L^{\frac{10}{3}}$ bounds,
\begin{align}
\label{twostars} \|\psi_j\|_{L^{\frac{10}{3}}_{x,t}} & \lesssim
N^{s-1} (N_j)^{-s} Z_I(T^*), \quad j \geq 1.
\end{align}

We now have the ingredients for our desired  $L^5_{x,t}$ bound of
$\phi$. By the triangle inequality,
\begin{align}
\label{sumthis} \|\phi\|_{L^5_{x,t}} & \leq  \sum_{j=0}^\infty
\|\psi_j\|_{L^5_{x,t}}.
\end{align}
Interpolating between the $L^{10}$ and $L^{4}$ bounds of
\eqref{star},\eqref{gotitgoinon} gives,
\begin{align}
\label{lowsdone} \|\psi_0\|_{L^5_{x,t}} & \lesssim \|\psi_0
\|^{\frac{2}{3}}_{L^4_{x,t}} \cdot
\| \psi_0\|_{L^{10}_{x,t}}^{\frac{1}{3}} \\
& \lesssim  \epsilon^{\frac{2}{3}} \left(Z_I(T^*)
\right)^{\frac{1}{3}}.
\end{align}
For $N_j \gtrsim   N$ interpolation between \eqref{star} and
\eqref{twostars} yields,
\begin{align*}
\sum_{j=1}^{\infty} \| \psi_j \|_{L^5_{x,t}} & \lesssim
\sum_{j=1}^\infty \|\psi_j
\|^{\frac{1}{2}}_{L^{\frac{10}{3}}_{x,t}}
\cdot \| \psi_j \|^{\frac{1}{2}}_{L^{10}_{x,t}} \\
& \lesssim \sum_{j=1}^\infty \left( N^{s-1} (N_j)^{-s} Z_I(T^*)
\right)^{\frac{1}{2}} \cdot
\left( (N_j)^{1-s} \cdot N^{s-1}  Z_I(T^*) \right)^{\frac{1}{2}} \\
& \lesssim N^{s-1}  Z_I(T^*),
\end{align*}
since $s > \frac{1}{2}$.  Choosing $N$ sufficiently large,
depending on $\epsilon$, yields \eqref{interp} for these high
frequency contributions as well. \qed
\medskip

{\noindent \bf Proof of Proposition \ref{almostconservation}}

For sufficiently smooth solutions, the usual energy \eqref{energy}
is shown to be conserved by differentiating in time, integrating
by parts, and using the equation \eqref{nls},
\begin{align*}
 \frac{d}{dt} E(\phi) &= \Re \int_{\rr^3} \overline{\phi_t} (|\phi|^2 \phi -
\Delta \phi) dx \\
& = \Re \int_{\rr^3} \overline{\phi_t} ( |\phi|^2 \phi -
\Delta \phi - i \phi_t) dx\\
& = 0.
\end{align*}
We begin to estimate $E(I\phi)(t)$ in  the same way. We need to
pay attention when we use the equation \eqref{nls} since of course
$I \phi$ is not a solution.  Repeating our steps above gives,
\begin{align*}
\frac{d}{dt} E(I\phi)(t) & = \Re \int_{\rr^3} \overline{I(\phi)_t}
( |I\phi|^2
I\phi - \Delta I\phi - iI\phi_t) dx \\
& = \Re \int_{\rr^3} \overline{I(\phi)_t} ( |I\phi|^2 I\phi -
I(|\phi|^2 \phi)) dx.
\end{align*}
When we integrate in time and apply the Parseval formula it
remains for us to bound
\begin{align}
E(I\phi(t))-E(I\phi(0))
  & =
\Re \int_0^t \int_{\sum_{j=1}^4 \xi_j = 0} \left ( 1 - \frac{
m(\xi_2 + \xi_3 + \xi_4)}{m(\xi_2) \cdot m(\xi_3) \cdot m(\xi_4)}
\right) \widehat{\overline{I \partial_t \phi}} (\xi_1)
\widehat{I\phi}(\xi_2) \widehat{\overline{I \phi}}(\xi_3)
\widehat{I \phi}(\xi_4). \label{fundII}
\end{align} We use the equation \eqref{nls} to substitute for
$\partial_t I(\phi)$ in \eqref{fundII}. Our aim is to show that
\begin{align}
\termone + \termtwo & \lesssim N^{- 1+} (Z_I(T))^P,
\label{oneandtwo}
\end{align}
for some $P > 0$, where the two terms on the left are
\begin{align}
\termone & \equiv \left |\int_0^T \int_{\sum_{i=1}^4 \xi_i = 0}
\left( 1 - \frac{m(\xi_2 + \xi_3 + \xi_4)}{m(\xi_2)  m(\xi_3)
m(\xi_4)} \right) \widehat{( \Delta \overline{I\phi})}(\xi_1)
\cdot \widehat{I\phi}(\xi_2) \cdot
\widehat{\overline{I\phi}}(\xi_3) \cdot \widehat{I\phi}(\xi_4)
\right|
\label{I} \\
\termtwo & \equiv \left| \int_0^T \int_{\sum_{i=1}^4 \xi_i = 0}
\left( 1 - \frac{m(\xi_2 + \xi_3 + \xi_4)}{m(\xi_2)  m(\xi_3)
m(\xi_4)} \right) \widehat{(\overline{I(|\phi|^2 \phi))}}(\xi_1)
\cdot \widehat{I\phi}(\xi_2) \cdot
\widehat{\overline{I\phi}}(\xi_3) \cdot \widehat{I\phi}(\xi_4)
\right|. \label{II}
\end{align}
In both cases we break $\phi$ into a sum of dyadic constituents
$\phi_j$ , each localized with a smooth cut-off function in
spatial frequency space to have support $\langle \xi \rangle
\thicksim 2^{k_j} \equiv N_j$, $k_j \in \{0, \ldots\}$ , and
employ the following estimate of Coifman-Meyer for a class of
multilinear operators.

Consider an infinitely differentiable symbol $\sigma: \rr^{nk}
\rightarrow \cc$ so that for all $\alpha \in N^{nk}$ and all $\xi
= (\xi_1 , \ldots , \xi_k) \in \rr^{nk}$, there is a constant
$c(\alpha)$ with,
\begin{align}
\label{symbolassume} |\partial^{\alpha}_{\xi} \sigma(\xi) | & \leq
c(\alpha) (1 + |\xi|)^{- |\alpha|}.
\end{align}
Define the multilinear operator $\Lambda$ by,
\begin{align}
\label{cmoperator} [\Lambda(f_1, \ldots, f_k)](x) & =
\int_{\rr^{nk}} e^{ix(\xi_1 + \ldots + \xi_k)} \sigma(\xi_1,
\ldots \xi_k) \hat f_1(\xi_1) \cdots \hat f_k(\xi_k) d \xi_1
\cdots d\xi_k.
\end{align}

\begin{theorem}[\cite{cmfourier}, Page 179]
\label{cmmultiplier} Suppose $p_j \in (1, \infty)$, $j = 1, \ldots
k$, are  such that $\frac{1}{p} = \frac{1}{p_1} + \frac{1}{p_2} +
\cdots + \frac{1}{p_k}  \, \leq 1$. Assume $\sigma(\xi_1, \ldots
\xi_k)$ a smooth symbol as in \eqref{symbolassume}.  Then there is
a constant $C = C(p_i,n,k,c(\alpha))$ so that for all Schwarz
class functions $f_1, \ldots f_k$,
\begin{align}
\label{multilinearestimate} \| \Lambda(f_1, \ldots, f_k)
\|_{L^p(\rr^n)} & \leq C \|f_1\|_{L^{p_1}(\rr^n)} \cdots
\|f_k\|_{L^{p_k}(\rr^n)}
\end{align}
\end{theorem}

\noindent Remark:  The estimate \eqref{multilinearestimate} is
also available for operators whose symbols obey much weaker bounds
than \eqref{symbolassume}, see e.g.  \cite{cmaster}, page 55.

When we   estimate below the terms  which constitute both
$\termone$ \eqref{I}  and $\termtwo$  \eqref{II}, we will first
seek a pointwise bound on the symbol,
\begin{align}
\label{symbol} \big| 1 -  \frac{m(\xi_2 + \xi_3 + \xi_4)}{m(\xi_2)
m(\xi_3) m(\xi_4)}  \big|& \leq B(N_2,N_3,N_4).
\end{align}
We factor $B(N_2,N_3,N_4)$ out of the left side of \eqref{symbol},
leaving a symbol $\sigma$ that satisfies the estimate
\eqref{symbolassume}\footnote{The required $L^\infty$ bound is
clear, and we leave the reader to check that the derivatives are
bounded as in \eqref{symbolassume}.}.  We are left to estimate  a
quantity of the form
\begin{equation*}
\left | B(N_2, N_3, N_4) \int_0^T \int_{\rr^3} [\Lambda(f_1,
f_2,f_3)]\text{\^{}}(\xi_4) \overline{\hat{f_4}}(\xi_4) d \xi_4 dt
\right|,
\end{equation*}
for some multilinear operator $\Lambda$ of the form
\eqref{cmoperator}, \eqref{symbolassume}. We estimate this using
the Plancherel formula, H\"older's inequality, Theorem
\ref{cmmultiplier}, and the Strichartz estimates. We can sum over
all the  dyadic pieces $\phi_j$ since our bounds will be seen to
decay sufficiently fast in the frequencies $N_i$.  We suggest that
the reader at first ignore this summation issue, and so ignore on
first reading the appearance below of all factors such as
$N_i^{0-}$ which we include only to show explicitly why our
frequency interaction estimates allow us to sum over the pieces
$\phi_i$. The main goal of the analysis is to establish the decay
of $N^{-1 +}$ in each class of frequency interactions below. In
what follows we drop the complex conjugates as they don't affect
the analysis used here{\footnote{A more detailed argument
exploiting the complex conjugates as in \cite{KPVQuad,
collianderfancy, taolove} might obtain a better exponent in
\eqref{increment}}}.

Consider first $\termone$.  We will conclude that $\termone \leq
N^{-1+}$ once we prove
\begin{multline} \label{termonegoal}
 \left| \int_0^T \int_{\sum_{i=1}^4 \xi_i = 0}
\left ( 1 - \frac{ m(\xi_2 + \xi_3 + \xi_4)}{m(\xi_2) \cdot
m(\xi_3) \cdot m(\xi_4)} \right) \widehat{\phi_1}(\xi_1)
\widehat{\phi_2}(\xi_2) \widehat{ \phi_3}(\xi_3)
\widehat{\phi_4}(\xi_4) \right|  \\
\lesssim \quad N^{- 1+} C(N_1, N_2, N_3, N_4) \left(Z_I(T)
\right)^4
\end{multline}
where  $C(N_1, N_2, N_3, N_4)$ is sufficiently small. By symmetry,
we may assume $N_2 \geq N_3 \geq N_4$.  The precise extent to
which $C(N_1, N_2, N_3, N_4)$ decays in its arguments, and the
fact that this decay allows us to sum over all dyadic shells, will
be described below.

\noindent{\bf $\termone$, Case 1:  $N \gg N_2$}.  According to
\eqref{Iproperties}, the symbol $1 - \frac{ m(\xi_2 + \xi_3 +
\xi_4)}{m(\xi_2) \cdot m(\xi_3) \cdot m(\xi_4)}$ on the right of
\eqref{fundII} is in this case identically zero and the bound
\eqref{termonegoal} holds trivially.

\noindent{\bf $\termone$, Case 2:  $N_2 \gtrsim N \gg N_3 \geq
N_4$}. Since $\sum_i \xi_i = 0$, we have $N_1 \thicksim N_2$. We
aim for \eqref{termonegoal} with
\begin{align}
\label{Ccase2_termone} C(N_1, N_2, N_3, N_4) & = N_2^{0-}.
\end{align}
With this decay factor, and the fact that we are considering here
terms where $N_1 \thicksim N_2$, we may immediately sum over all
the $N_i$.

By the mean value theorem,
\begin{align}
\left| \frac{m(\xi_2) - m(\xi_2 + \xi_3 + \xi_4)}{m(\xi_2)}
\right| & \lesssim \frac{|\nabla m(\xi_2)\cdot(\xi_3 +
\xi_4)|}{m(\xi_2)} \lesssim \frac{N_3}{N_2} .\label{pointwiseII}
\end{align}
After estimating the symbol with \eqref{pointwiseII}, we view the
$N_3$ in the numerator as resulting from a derivative falling on
the $I \phi_3$ factor in the integrand.  Hence these interactions
can be estimated using  H\"older's inequality, Theorem
\ref{cmmultiplier},  and the definition \eqref{newz} of $Z_I(t)$,
\begin{align*}
\left|\hbox{Left Side  of \eqref{termonegoal}} \right| & \lesssim
\frac{N_3}{N_2} \left| \int_0^T \int_{\rr^3}
\Lambda[\Delta I \phi_1, I \phi_2 , I \phi_3 ] \cdot I \phi_4 dx dt \right|\\
& \leq \frac{1}{N_2} || \Delta I \phi_1
||_{L^{\frac{10}{3}}_{x,t}} \cdot || I \phi_2
||_{L^{\frac{10}{3}}_{x,t}}
\cdot ||\nabla I \phi_3||_{L^{\frac{10}{3}}_{x,t}} \cdot ||I \phi_4||_{L^{10}_{x,t}} \\
& \leq \frac{N_1}{N_2 \cdot N_2} \cdot (Z_I(t))^4 \\
& \leq \frac{1}{N_1} (Z_I(t))^4 \\
& \leq N^{-1+} \cdot N_2^{0-} (Z_I(t))^4
\end{align*}
by our assumptions on the $N_i$.  This establishes
\eqref{termonegoal}, \eqref{Ccase2_termone}.

\noindent{\bf $\termone$, Case 3:  $N_2 \geq N_3 \gtrsim N $}.  In
this case the only
pointwise bound available for the symbol is the straightforward
one:  when $|\xi_1|, |\xi_2|$ are not comparable, no cancellation
can occur in the numerator of \eqref{symbol}. When $|\xi_1|
\thicksim |\xi_2|$, we then also need $|\xi_3|, |\xi_4| \leq N$ in
order to get cancellation.  If any of these conditions fail, our
pointwise estimate will be simply,
\begin{align}
\label{dumbbound} \left| 1 - \frac{m(\xi_2 + \xi_3 + \xi_4)}{
m(\xi_2) m(\xi_3) m(\xi_4)} \right| &
\lesssim\frac{m(\xi_1)}{m(\xi_2)m(\xi_3) m(\xi_4)}.
\end{align}
The frequency interactions here fall into two subcategories,
depending on which frequency is comparable to $N_2$. {\noindent
{\bf Case 3(a):} $N_1 \sim N_2 \geq N_3 \gtrsim N$.} By
assumption, $s > \frac{1}{2} + \delta$ for some small $\delta$. In
this case we prove the decay factor
\begin{align}
\label{Ccase3a_termone} C(N_1,N_2,N_3,N_4) & = N^{-1+2 \delta}
N_3^{0- 2 \delta}
\end{align}
in \eqref{termonegoal}.  This allows us to directly sum in
$N_3,N_4$, and sum in $N_1, N_2$ after applying Cauchy-Schwarz to
those factors. Estimate the symbol using \eqref{dumbbound}. Use
H\"older's inequality and Theorem \ref{cmmultiplier} to take the
factors involving $\phi_i, i = 1,2,3$ in $L^{\frac{10}{3}}_{x,t}$,
and the $\phi_4$ factor in $L^{10}_{x,t}$.  It remains to show
\begin{align}
\label{case3astuff} \frac{m(N_1) N_1 N^{1-2 \delta} N_3^{2
\delta}}{m(N_2) m(N_3)m(N_4)N_2N_3} & \lesssim 1.
\end{align}
When proving such estimates here and in the sequel, we shall
frequently use the following two elementary facts without further
mention:  for any $p > \frac{1}{2} - \delta$, the function $m(x)
|x|^p$ is increasing, and $m(x) \langle x \rangle$ is bounded
below.  The bound \eqref{case3astuff} is now straightforward,
\begin{align*}
{\text{Left Side of \eqref{case3astuff}}} & \lesssim \frac{N^{1-2 \delta} N_3^{2 \delta}}{m(N_3)m(N_4)N_3} \\
& \lesssim \frac{N^{1- 2 \delta} N_3^{2 \delta}}{(m(N_3))^2 N_3} \\
& \lesssim \frac{N^{1- 2 \delta} N_3^{2 \delta}}{(m(N_3))N_3^{\frac{1}{2} - \delta} m(N_3) N_3^{\frac{1}{2} - \delta} N_3^{2 \delta}}  \\
& \lesssim \frac{N^{1-2 \delta} N_3^{2 \delta}}{N^{1-2 \delta}
N_3^{2 \delta}}
\end{align*}
which gives \eqref{termonegoal}, \eqref{Ccase3a_termone}.

{\noindent {\bf Case 3(b):} $N_2 \sim N_3  \gtrsim N$.} We aim in
this case for the decay factor
\begin{align}
\label{Ccase3b_termone} C(N_1,N_2,N_3,N_4) & = N^{-1+ 2 \delta }
N_2^{- 2 \delta}
\end{align}
where $\delta$ is as in Case 3(a) above.  This will allow us to
sum directly in all the $N_i$.  Once again we use
\eqref{dumbbound} and apply H\"older's inequality and
\eqref{multilinearestimate} exactly as in the preceding
discussion.
\begin{align*}
\frac{m(N_1) N_1 N^{1-2 \delta} N_2^{2 \delta}}{m(N_2)
m(N_3)m(N_4)N_2N_3} & \lesssim
\frac{m(N_1) N_1 N^{1- 2 \delta} N_2^{2 \delta}}{(m(N_2))^3 N_2N_2}\\
& \lesssim \frac{m(N_2) N_2 N^{1-2 \delta} N_2^{2 \delta}}{(m(N_2))^3 N_2N_2} \\
& = \frac{N^{1- 2 \delta} N_2^{2 \delta}}{(m(N_2))^2\cdot N_2} \\
& \leq \frac{N^{1-2 \delta} N_2^{2 \delta}}{N_2^{2 \delta} \cdot N^{1-2 \delta}} \\
& \leq 1,
\end{align*}
as desired.  It remains to prove bounds of the form
\eqref{oneandtwo} for $\termtwo$\eqref{II}.

When decomposing the integrand of $\termtwo$ in frequency space,
write $N_{123}$ for the dyadic frequency into which we project the
nonlinear factor $I(\phi^3)$. Note that in the treatment of
$\termone$ above, we always took the $\Delta \phi_1$ factor in
$L^{\frac{10}{3}}$, estimating this by $N_1 Z_I(T)$. The analysis
above for $\termone$ therefore applies unmodified to $\termtwo$
once we prove the following,
\begin{lemma} \label{termtwolemma}  Assume $\phi, T, Z_I(T), N_{123}$ as defined above, and $P_{N_{123}}$
the Littlewood-Paley projection onto the $N_{123}$ frequency
shell.  Then
\begin{align}
\label{termtwostar} \| P_{N_{123}} (I(\phi^3))
\|_{L^{\frac{10}{3}}_{x,t}([0,T]\times \rr^3)} & \lesssim N_{123}
(Z_I(T))^3.
\end{align}
\end{lemma}
\noindent{{\bf Proof:}}  We write $\phi = \phi_L + \phi_H$ where
\begin{align*}
\text{supp} \hat{\phi}_L(\xi,t) & \subseteq \{ |\xi| < 2 \} \\
\text{supp} \hat{\phi}_H(\xi,t) & \subseteq \{ |\xi| > 1 \}.
\end{align*}
Consider first the bound \eqref{termtwostar} when all three
factors on the left are $\phi_L$,
\begin{align*}
\| P_{N_{123}} (I(\phi_L^3)) \|_{L^{\frac{10}{3}}_{x,t}} & \lesssim \| \phi_L \|_{L^{10}_{x,t}}^3 \\
& = \| I \phi_L \|_{L^{10}_{x,t}}^3 \\
& \leq (Z_I(T))^3 \\
& \lesssim N_{123} (Z_I(T))^3,
\end{align*}
since  $N_{123} \geq 1$. When instead all three components on the
left of \eqref{termtwostar} are $\phi_H$,  we have by
Littlewood-Paley theory, Sobolev embedding, and the Leibniz rule
mentioned in the proof of Proposition \ref{almostconservation},
\begin{multline*}
\| \frac{1}{N_{123}} P_{N_{123}} I(\phi_H^3)
\|_{L^{\frac{10}{3}}_{x,t}} \; \lesssim \;  \| \nabla^{-1}
P_{N_{123}}
I(\phi_H^3) \|_{L^{\frac{10}{3}}_{x,t}} \\
\lesssim \: \| \nabla^{\frac{1}{2}}  I (\phi_H^3)
\|_{L^{\frac{10}{3}}_t L^{\frac{10}{8}}_x}
\; \lesssim \;  \|\nabla^{\frac{1}{2}} I \phi_H \|^3_{L^{10}_tL^{\frac{30}{8}}_x} \\
\; \lesssim \; \| \nabla I \phi_H \|^3_{L^{10}_t
L^{\frac{30}{13}}_x} \; \lesssim \; (Z_I(T))^3 \quad \quad
\end{multline*}
as desired.

The remaining terms are bounded using similar arguments,
\begin{multline*}
 \|\frac{1}{N_{123}} P_{N_{123}} I(\phi_H \cdot \phi_H \cdot \phi_L) \|_{L^{\frac{10}{3}}_{x,t}}  \; \;
\lesssim  \; \;  \|\nabla^{\frac{1}{2}} I ( \phi_H \cdot \phi_H \cdot \phi_L) \|_{L^{\frac{10}{3}}_t L^{\frac{10}{8}}_x} \\
\;\lesssim \; \| \nabla^{\frac{1}{2}} I
\phi_H\|_{L^{10}_tL^{\frac{30}{8}}_x} \cdot
\|\phi_H\|_{L^{10}_tL^{\frac{30}{13}}_x} \cdot \| \phi_L
\|_{L^{10}_tL^{10}_x}  + \| \phi_H \|_{L^{10}_t
L^{\frac{30}{8}}_x} \cdot \| \phi_H \|_{L^{10}_t
L^{\frac{30}{8}}_x} \cdot
\| \nabla^{\frac{1}{2}} I  \phi_L \|_{L^{{10}_t L^{\frac{30}{8}}_x}} \\
\;\lesssim \; \|\nabla I \phi_H \|_{L^{10}_tL^{\frac{30}{13}}_x}
\cdot \|\nabla I \phi_H\|_{L^{10}_tL^{\frac{30}{13}}_x} \cdot \| I
\phi_L  \|_{L^{10}_tL^{10}_x} + \| \nabla^{\frac{1}{2}} I \phi_H
\|_{L^{10}_t L^{\frac{30}{8}}_x} \cdot \| \nabla^{\frac{1}{2}} I
\phi_H \|_{L^{10}_t L^{\frac{30}{8}}_x}\cdot
\| \nabla I  \phi_L \|_{L^{10}_t L^{\frac{30}{13}}_x}\\
\;  \lesssim \; (Z_I(T))^3. \quad \quad\quad \quad\quad \quad\quad
\quad
 \end{multline*}
 \begin{multline*}
  \|\frac{1}{N_{123}} P_{N_{123}} I(\phi_H \cdot \phi_L \cdot \phi_L) \|_{L^{\frac{10}{3}}_{x,t}}
  \; \lesssim \; \| \phi_H \cdot \phi_L \cdot \phi_L \|_{L^{\frac{10}{3}}_t L^{\frac{30}{19}}_x}  \\
  \lesssim \|\phi_H\|_{L^{10}_tL^{\frac{30}{13}}_x} \cdot \|\phi_L \|_{L^{10}_t L^{10}_x}
  \cdot \| \phi_L\|_{L^{10}_t L^{10}_x} \; \lesssim \; \|\nabla I \phi_H\|_{L^{10}_tL^{\frac{30}{13}}_x}
  \cdot  \| \nabla I \phi_L\|_{L^{10}_t L^{\frac{30}{13}}_x}^2 \\
  \lesssim \; (Z_I(T))^3. \quad \quad \quad \quad
\end{multline*}
This completes the proof of Lemma \ref{termtwolemma}, and hence
Proposition \ref{almostconservation}. \qed

\section{Proof of Main Theorem}

We combine the interaction Morawetz estimate
\eqref{spacetimeLfour} and Proposition \ref{almostconservation}
with a scaling argument to prove the following statement giving
uniform bounds in terms of the rough norm of the initial data.

\begin{proposition} \label{globalL4bound}
Suppose $\phi(x,t)$ is a global in time solution to
\eqref{nls}-\eqref{nlsdata} from data $\phi_0 \in
C^{\infty}_0(\rr^3)$.  Then so long as $s > \frac{4}{5}$, we have
\begin{align}
\label{reducedtoshow}
||\phi||_{L^4([0,\infty] \times \rr^3)} & \leq C(||\phi_0||_{H^s(\rr^3)}) \\
\sup_{0 \leq t < \infty} ||\phi(t)||_{H^s(\rr^3)} & \leq
C(||\phi_0||_{H^s(\rr^3)}). \label{uniformHs}
\end{align}
\end{proposition}


\noindent{\bf Remark:}  As mentioned at the outset of the paper,
energy conservation \eqref{energy} and the local in time
well-posedness of \eqref{nls}-\eqref{nlsdata} from data in
$H^s(\rr^3)$, $s > \frac{1}{2}$  imply that the solution $\phi$
considered here is smooth and exists globally in time.  Since the
estimate \eqref{uniformHs} involves only the rough norm
$||\phi_0||_{H^s(\rr^3)}$ on the right hand side, the global
well-posedness portion of Theorem \ref{maintheorem} follows from
\eqref{uniformHs}, the local existence theory (see \cite{cw:local}
for a proof and further references), and a standard density
argument.

\begin{proof}

The first step is to scale the solution:  if $\phi$ is a solution
to \eqref{nls}, then so is
\begin{align}
\label{scaling} \philambda(x,t) & \equiv \frac{1}{\lambda}
\phi(\frac{x}{\lambda}, \frac{t}{\lambda^2}).
\end{align}
We choose $\lambda$ so that
 $E(I \philambda_0) \, \equiv \, \frac{1}{2} ||\nabla I \philambda_0 ||^2_{L^2(\rr^3)} +
\frac{1}{4} \| I \philambda_0 \|_{L^4_x}^{4} \, \leq \,
\frac{1}{4}$. This is possible since we are working with
subcritical $s$, so long as we choose $\lambda$ in terms of the
parameter\footnote{The parameter $1 \ll N $ will be chosen at the
very end of the argument, where it is shown to depend only on
$||\phi_0||_{H^s(\rr^3)}$ .} $N$.  Specifically, arguing as in
\eqref{boundenergy}, one easily shows,
\begin{align*}
\frac{1}{2} ||\nabla I \philambda_0 ||^2_{L^2(\rr^3)} & \lesssim
\left( N^{1 - s} \lambda ^{\frac{1}{2} - s}
\|\phi_0\|_{H^s(\rr^3)} \right)^2
\end{align*}
In order to make the right hand side here $\leq \frac{1}{8}$,
choose
\begin{align}
\label{chooselambda} \lambda \approx
N^{\frac{1-s}{s-\frac{1}{2}}}.
\end{align}
One can bound the second term in $E(I \philambda_0)$ by
considering separately the domains $|\xi| \lesssim
\frac{1}{\lambda}, \frac{1}{\lambda} \lesssim |\xi| \lesssim N$,
and $ |\xi| \gtrsim N$ in frequency space: straightforward
arguments using Sobolev embedding together with the relation
\eqref{chooselambda} will give
\begin{align*}
\frac{1}{4} \| I \philambda_0 \|_{L^4_x}^{4} & \leq \frac{1}{8}.
\end{align*}

 We claim that the set $W$ of times for which
\eqref{reducedtoshow} holds is all of $[0,\infty)$.  In the
process of proving this, we will also show \eqref{uniformHs} holds
on $W$.

For some universal constant $C_1$ to be chosen shortly,
define\footnote{Roughly speaking, our bound for $||
\philambda||_{L^4([0,T] \times \rr^3)}$ in this definition scales
like $\lambda^{\frac{3}{8}}$ as the $L^4_{x,t}$ estimate provide
by \eqref{spacetimeLfour} has - just looking at low frequency
contributions for the moment- $\frac{3}{4}$ of a factor of $\|
P_{\leq N} \philambda \|_{L^2_x}$ and $\frac{1}{4}$ of a factor of
$\|\nabla P_{\leq N} \philambda \|_{L^2_x}$ on the right hand
side.  The former scales like
$(\lambda^{\frac{1}{2}})^{\frac{3}{4}}$, while a bootstrap
argument will show the latter is $\leq 1$.}
\begin{align}
\label{l4claim} W & \equiv \left\{ T \; : \; ||
\philambda||_{L^4([0,T] \times \rr^3)}\; \leq \; C_1
\lambda^{\frac{3}{8}} \right\} .
\end{align}
The set $W$ is clearly closed and nonempty. It suffices then to
show it is open. Note that the  quantity
$\|\philambda\|_{L^4([0,T] \times \rr^3)}$ is continuous in time
as we've reduced to the case when $\phi(x,t)$ is smooth.  Hence if
$T_1 \in W$, then for some $T_0 > T_1$ sufficiently close to $T_1$
we have
\begin{align}
\label{continuity} ||\philambda||_{L^4_{x,t}([0,T_0] \times
\rr^3)} & \leq 2 C_1 \lambda^{\frac{3}{8}}.
\end{align}
We claim $T_0 \in W$.  By \eqref{spacetimeLfour},
\begin{align}
|| \philambda||_{L^4_{x,t}([0,T_0] \times \rr^3)} & \lesssim
||\philambda_0||^{\frac{1}{2}}_{L^2_x} \cdot \sup_{0 \leq t \leq
T_0}||\philambda
(t)||_{\dot{H}^{\frac{1}{2}}(\rr^3)}^{\frac{1}{2}}.
\label{firststep} \\
& \leq  C(\|\phi_0\|_{L^2_x}) \lambda^{\frac{1}{4}} \cdot \sup_{0
\leq t \leq T_0} || \philambda (t) ||^{\frac{1}{2}}_{
\dot{H}^{\frac{1}{2}}(\rr^3)} \label{secondstep}
\end{align}
where we've taken into account the $L^2$ conservation law
\eqref{l2conservation}.  To bound the second factor in
\eqref{secondstep}, decompose $\philambda (t)$ as,
\begin{align}
\label{breakitupbreakitup} \philambda(t) & = P_{\leq N} \philambda
(t) + P_{\geq N} \philambda (t).
\end{align}
That is, a sum of functions supported on frequencies $|\xi| \leq
N$ and $|\xi| \geq N$, respectively.  Interpolation and the fact
that $I$ is the identity on low frequencies gives us the bound,
\begin{align}
\|P_{\leq N} \philambda (t) \|_{\dot{H}^{\frac{1}{2}}_x} &
\lesssim \|P_{\leq N} \philambda (t)\|_{L^2_x}^{\frac{1}{2}} \cdot
\| P_{\leq N} \philambda (t)\|_{\dot{H}^1_x}^{\frac{1}{2}} \nonumber \\
& \lesssim  \| \philambda_0 \|_{L^2_x}^{\frac{1}{2}} \cdot
\| I P_{\leq N} \philambda (t)\|_{\dot{H}^1_x}^{\frac{1}{2}} \nonumber \\
& \leq C(\|\phi_0 \|_{L^2_x}) \lambda^{\frac{1}{4}} \| I
\philambda (t) \|_{\dot{H}^1_x}^{\frac{1}{2}}. \label{gotI}
\end{align}
We interpolate the high frequency constituent between
$\dot{H}^s_x$ and $L^2_x$, and use the definition
\eqref{Iproperties} of $I$ to get,
\begin{align}
\|P_{\geq N} \philambda (t) \|_{\dot{H}^{\frac{1}{2}}_x} &
\lesssim \|P_{\geq N} \philambda (t)\|_{L^2_x}^{1 - \frac{1}{2s}}
\cdot
\| P_{\geq N} \philambda (t)\|_{\dot{H}^s_x}^{\frac{1}{2s}} \nonumber \\
& = \|P_{\geq N} \philambda (t)\|_{L^2_x}^{1 - \frac{1}{2s}} \cdot
N^{\frac{s-1}{2s}}
\| I P_{\geq N} \philambda (t)\|_{\dot{H}^1_x}^{\frac{1}{2s}} \nonumber \\
& \leq C(\|\phi_0\|_{L^2_x}) \cdot \|I \philambda
\|_{\dot{H}^1_x}^{\frac{1}{2s}}, \label{gotII}
\end{align}
where we've used both the $L^2$ conservation
\eqref{l2conservation} and our choice of $\lambda$,
\eqref{chooselambda}. Putting together \eqref{gotII},
\eqref{gotI}, \eqref{breakitupbreakitup}, and \eqref{secondstep}
gives us
\begin{align}
\label{thisiswherewearelikeitornot} ||
\philambda||_{L^4_{x,t}([0,T] \times \rr^3)}   & \leq C(\| \phi_0
\|_{L^2_x}) \bigg( \lambda^{\frac{3}{8}}
 \sup_{0 \leq t \leq T_0} || I \philambda (t) ||^{\frac{1}{4}}_{\dot{H}^{1}_x} +
 \sup_{0 \leq t \leq T_0} || I \philambda (t) ||^{\frac{1}{4s}}_{\dot{H}^{1}_x} \bigg).
\end{align}

We conclude $T_0 \in W$ if we establish
\begin{align}
\label{toshowII} \sup_{0 \leq t \leq T_0} \| I \philambda (t) \|_{
\dot{H}^1(\rr^3)} & \leq 1
\end{align}
since we then take $C_1$ in \eqref{l4claim} larger than twice the
constant $C(\| \phi_0 \|_{L^2_x})$ appearing in
\eqref{thisiswherewearelikeitornot}.

By \eqref{continuity} we may divide the time interval $[0,T_0]$
into subintervals $I_j, j = 1, 2, \ldots,L$ so that for each $j$,
\begin{align}
\label{cutsmall} ||\philambda||_{L^4_{x,t}(I_j \times \rr^3)} &
\leq \epsilon.
\end{align}
Apply the almost conservation law in Proposition
\ref{almostconservation} on each of the subintervals $I_j$ to get
\begin{align}
\label{lotsofbumps} \sup_{0 \leq t \leq T_0} ||\nabla I
\philambda(t) ||_{L^2(\rr^3)} & \leq  E(I\phi_0) + C L \cdot
N^{-1+}.
\end{align}
We get \eqref{toshowII} from \eqref{lotsofbumps} if we can show
\begin{align}
L \cdot N^{-1+} & \ll \frac{1}{4}. \label{needy}
\end{align}
Recall $L$ was defined essentially by \eqref{cutsmall}.  Since
\begin{align*}
||\philambda||^4_{L^4_{x,t}([0,T_0] \times \rr^3)} & \lesssim
\lambda^{\frac{3}{2}},
\end{align*}
we can be certain that $L \approx \lambda^{\frac{3}{2}}$.  If we
put this together with \eqref{needy} and \eqref{chooselambda}, we
see that we need to be able to choose $N$ so that
\begin{align*}
(N^{\frac{1-s}{s-\frac{1}{2}}})^{\frac{3}{2}} \cdot N^{-1+} & \ll
\frac{1}{4}.
\end{align*}
This is possible since  for $s > \frac{4}{5}$ the exponent on the
left is negative. Notice that \eqref{uniformHs} holds on the set
$W$ using \eqref{toshowII}, the definition of $I$, and $L^2$
conservation.

\end{proof}

We have already explained  why the global well-posedness statement
in Theorem \ref{maintheorem} follows from \eqref{uniformHs}.  It
remains only to prove scattering using the following well-known
arguments.  (See  e.g. \cite{linstrauss, GVScatter, BScatter,
cazbook}.) Asymptotic completeness will follow quickly once we
establish a uniform bound of the form,
\begin{align}
\label{z}
Z(t) & \equiv \sup_{q,r \; \hbox{admissible}} || \langle \nabla \rangle^s \phi ||_{L^q_tL^r_x([0,t] \times \rr^3)} \\
& \leq C(||\phi_0||_{H^s(\rr^3)}). \label{boundz}
\end{align}
This is established much as in the proof of Lemma
\ref{interpolation}. By \eqref{reducedtoshow}, we can  decompose
the time interval $[0, \infty)$ into a finite number of disjoint
intervals $J_1, J_2, \ldots J_K$ where for $i = 1, \ldots K$ we
have
\begin{align}
\label{smalll4norm} ||\phi||_{L^4_{x,t}(J_i \times \rr^3)} & \leq
\epsilon
\end{align}
for a constant $\epsilon(||\phi_0||_{H^s(\rr^3)})$ to be chosen
momentarily.

Apply $\langle \nabla \rangle ^s$ to both sides of \eqref{nls}.
Choosing $\tilde{q}', \tilde{r}' = \frac{10}{7}$, the Strichartz
estimates \eqref{strichartzestimate} give us that for all $t \in
J_1$,
\begin{align*}
Z(t) & \lesssim ||\langle \nabla\rangle^s  \phi_0||_{L^2(\rr^3)} +
||\langle \nabla \rangle ^s  (\phi \phi
\phi)||_{L^{\frac{10}{7}}_{t,x}([0,t] \times \rr^3)}.
\end{align*}
Apply the fractional Leibniz rule to the last term on the right,
taking the factor with $\langle \nabla \rangle^s$ in
$L^{\frac{10}{3}}$, and the other two in  $L^5$.  The factor
ending up in $L^{\frac{10}{3}}$ is bounded by $Z(t)$.  The
remaining $L^5_{x,t}$ factors are bounded by interpolating between
$||\phi||_{L^4_{x,t}}$ and $||\phi||_{L^6_{x,t}}$.  The latter
norm is bounded by $Z(t)$ using Sobolev embedding:
$$||\phi||_{L^6_{x,t}} \; \lesssim \;  ||\langle \nabla \rangle^{\frac{2}{3}} \phi||_{L^6_t L^{\frac{18}{7}}_x} \;
\leq \; Z(t).$$ We conclude
\begin{align}
Z(t) & \lesssim  ||\phi_0||_{H^s(\rr^3)} + \epsilon^{\delta_1}
Z(t)^{(1 + \delta_2)}. \label{almostclosed}
\end{align}
for some constants $\delta_1, \delta_2 > 0$.  For sufficiently
small choice of $\epsilon$, the bound \eqref{almostclosed} yields
\eqref{boundz} for all $t \in J_1$, as desired.  Since we are
assuming the bound \eqref{uniformHs}, we may repeat this argument
to handle the remaining intervals $J_i$.

The asymptotic completeness claim  in Theorem \ref{maintheorem}
follows quickly from \eqref{boundz}. Given $\phi_0 \in
H^s(\rr^3)$, we look for a $\phi^+$ satisfying
\eqref{asymptoticstates}.  Set,
\begin{align}
\label{assigndata} \phi^+ & \equiv \phi_0 - i \int_0^{\infty}
S^L(-\tau) \left(|\phi|^2\phi \right) d\tau
\end{align}
\medskip
which will make sense once we show the integral on the right hand
side converges in $H^s(\rr^3)$.  Equivalently, we want
\begin{align}
\lim_{t\rightarrow \infty} || \int_t^{\infty} \langle \nabla
\rangle ^s S^L(-\tau) \left( |\phi|^2 \phi \right) d \tau
||_{L^2(\rr^3)} & = 0. \label{timeintegralbound}
\end{align}
With this,
\begin{align*}
\lim_{t \rightarrow \infty} ||S^L(t) \phi^+ -
\phi(t)||_{H^s(\rr^3)} & = \lim_{t \rightarrow \infty} ||\langle
\nabla \rangle^{s}
S^{L}(t) \int_t^{\infty} S^L(-\tau) \left( |\phi|^2\phi \right) d \tau ||_{L^2(\rr^3)} \\
& = 0
\end{align*}
since we are assuming \eqref{timeintegralbound}.    To prove
\eqref{timeintegralbound}, test the time integral on the left
against an arbitrary $L^2(\rr^3)$ function $F(x),
\|F(x)\|_{L^2(\rr^3)} \leq 1$. Using the fractional Leibniz rule,
\begin{align*}
\sup_{\|F(x)\|_{L^2(\rr^3)} \leq 1} \left\langle F(x) \; , \;
\int_t^{\infty} \langle \nabla \rangle^s S^L(-\tau) \left(
|\phi|^2 \phi \right) d \tau \right\rangle_{L^2(\rr^3)} & \approx
\sup_{\|F(x)\|_{L^2(\rr^3)} \leq 1} \left \langle S^L(\tau) F(x)
\; , \; (\nabla^s \phi) \phi \phi \right\rangle_{L^2_{x,t}([t,
\infty) \times
\rr^3)} \\
& \leq \sup_{\|F(x)\|_{L^2(\rr^3)} \leq 1} ||S^L(\tau)
F(x)||_{L^{\frac{10}{3}}_{x,\tau}} ||\nabla^s
\phi||_{L^{\frac{10}{3}}_{x,t}} ||\phi||^2_{L^5_{x,t}([t, \infty)
\times
\rr^3)}  \\
& \rightarrow 0,
\end{align*}
where in the last step the convergence is uniform in $F$, and
where we've used \eqref{boundz} and the $L^5_{x,t}$ argument
before \eqref{almostclosed}.  The statement
\eqref{timeintegralbound} follows by the converse to H\"older's
inequality.

For completeness we include an argument proving the existence of
wave operators on $H^s(\rr^3)$, following closely the exposition
of \cite{GVScatter} in  \cite{cazbook}, \S 7.6.  Given $\phi^+ \in
H^s(\rr^3)$, we are looking for a solution $\phi(x,t)$ of
\eqref{nls} and data $\phi_0$ which, heuristically at least,
satisfy,
\begin{align}
\phi(x,t) & = S^L(t) \phi_0 - i \int_0^t S^L(t-\tau) |\phi|^2 \phi d \tau \\
& = S^L(t) \left( S^{NL}(-\infty) S^L (\infty) \phi^+ \right) - i \int_0^t S^L(t-\tau) |\phi|^2 \phi d \tau \nonumber \\
& = S^L(t) \left( \phi^+ - i \int_{\infty}^0 S^L(0 - \tau)
|\phi|^2 \phi d \tau  \right) -
i \int_0^t S^L(t-\tau) |\phi|^2 \phi d \tau \nonumber \\
& = S^L(t) \phi^+ + i \int_t^{\infty} S^L(t - \tau) |\phi|^2 \phi
d \tau. \label{woequation}
\end{align}
 Heuristics aside, we now sketch how this last integral equation is solved for $\phi(x,t)$
using a fixed point argument, and prove that $\phi(x,t)$  does in
fact approach $S^L(t) \phi^+$ as $t \rightarrow \infty$.


By Strichartz estimates, we have $ S^L(t) \phi^+ \in
L^{\frac{8}{3}}_t W^{s,4}_x \cap L^{8}_t W^{s,\frac{12}{5}}_x ([0,
\infty) \times \rr^3)$.  Set,
\begin{align}
\label{ks} K_{t_0} & = \| S^L(t) \phi^+
\|_{L^{\frac{8}{3}}_tW^{s,4}_x([t_0, \infty) \times \rr^3)} + \|
S^L(t) \phi^+ \|_{L^{8}_tW^{s,\frac{12}{5}}_x([t_0, \infty) \times
\rr^3)}.
\end{align}
Clearly  $K_{t_0} \rightarrow 0$ as $t_0 \rightarrow \infty$.
Define,
\begin{align}
X & = \left\{ u \in L^{\frac{8}{3}}_tW^{s,4}_x \cap  L^{8}_t
W^{s,\frac{12}{5}}_x ((t_0, \infty) \times \rr^3 ) \; | \; \|
u\|_{L^{\frac{8}{3}}_t W^{s,4}_x((t_0, \infty) \times \rr^3)} + \|
u\|_{L^{8}_t W^{s,\frac{12}{5}}_x((t_0, \infty) \times \rr^3)} \;
\leq 2 K_{t_0} \right\} \label{X}
\end{align}
with norm $\| \cdot \|_{L^{\frac{8}{3}}_t W^{s,4}_x} +\| \cdot
\|_{L^{8}_t W^{s,\frac{12}{5}}_x} $.  For functions $u \in X$ we
have
\begin{align}
\| |u|^2 u \|_{L^{\frac{8}{5}}_t W^{s,\frac{4}{3}}_x} & \leq \|
\nabla^s u\|_{L^{\frac{8}{3}}_tL^4_x}
\cdot \|u\|^2_{L^8_tL^4_x} \label{algebraicboundI} \\
& \leq C (2 K_{t_0})^3 \label{algebraicboundII},
\end{align}
where we've bounded the second two factors on the right of
\eqref{algebraicboundI} using Sobolev embedding. It is
straightforward\footnote{The proof of Corollary 3.2.7 in
\cite{cazbook} can be followed without modification.} to conclude
from \eqref{algebraicboundII} that the function
\begin{align}
\label{forcedguy} \Phi_u(t) & \equiv i \int_t^{\infty} S^L(t -
\tau) |u|^2u d\tau
\end{align}
is well defined for all $u \in X$, and that
\begin{align}
\Phi_u(t) & \in C \left( (t_0,\infty); H^s(\rr^3) \right) \cap X,
\end{align}
with,
\begin{align}
\|\Phi_u \|_{X} & \leq C(2K_{t_0})^3 \leq K_{t_0},
\label{fourstars}
\end{align}
when  $K_{t_0}$ is small enough - that is, for $t_0$ large enough.
Hence the map,
\begin{align}
A: u(t) & \rightarrow S^L(t) \phi^+ + \Phi_u(t),
\end{align}
takes $X$ into itself.  It can be similarly argued that $A$ is a
contraction.  We conclude there is a unique solution $\phi \in X$
of \eqref{woequation}.  By our global existence result and time
reversibility, we may extend this solution $\phi$, starting from
data at time $t_0$, to all of $[0, \infty)$.  It is now
straightforward to verify that
$$
\lim_{t \rightarrow \infty} \| \phi(t) -
S^L(t)\phi^+\|_{H^s(\rr^3)} \; = \; 0,
$$
as desired. \qed


\appendix                               









\frenchspacing
\bibliographystyle{plain}

\end{document}